\documentclass{article} %{ndjflart} 

\usepackage{a4wide}
\usepackage{url}
\usepackage{amsmath,amssymb}
\usepackage[numbers]{natbib}
\usepackage{graphicx}
\usepackage[frame,curve]{xypic}
\usepackage{accents}
\newcommand{\sluttSymbol}{\Box} %{\blacksquare}
\newcommand{\slutt}{\hspace*{\fill $\sluttSymbol$}\\[1ex]}

\newenvironment{Proof}{\vspace*{.5ex}\par\noindent{\sc{Proof.}}}{\slutt}
\newenvironment{Proof*}[1]{\vspace*{.5ex}\par\noindent{\sc{Proof #1.}}}{\slutt}
\newtheorem{Theorem}{Theorem}[section]
\newtheorem{Lemma}[Theorem]{Lemma}
\newtheorem{Corollary}[Theorem]{Corollary}
\newtheorem{Fact}[Theorem]{Fact}

\newtheorem{Example}[Theorem]{Example}

\newcommand{\eq}[1]{\vspace*{-.5ex}\refstepcounter{Theorem}\begin{equation}#1\vspace*{-.5ex}\end{equation}}
\newcommand{\refp}[1]{(\ref{#1})}

\newcommand{\draw}[1]{\raisebox{.7ex}{\xymatrix@C=.35cm{#1}}}
\newcommand{\drawr}[2]{\raisebox{.7ex}{\xymatrix@C=.35cm@R=#1{#2}}}
\newcommand{\drawx}[2]{\raisebox{1}{\xymatrix@C=.35cm{#2}}}
\newcommand{\okk}[1]{#1^{ok}}

\newcounter{DELPOI}[Theorem]

\newcounter{AUX}

\newcommand{\ovr}[1]{\overline{#1}}
\newcommand{\any}[1]{\ddot{#1}} %{\tilde{#1}} % possibly: {\widetilde{#1}}
\newcommand{\ecl}{\{\}} %empty clause

\newcommand{\A}[1]{A_{\scriptscriptstyle{#1}}}

\newcommand{\Ar}[1]{A^*_{\scriptscriptstyle{#1}}}
\newcommand{\invA}[1]{A^-_{\scriptscriptstyle{#1}}}
\newcommand{\invrA}[1]{A^{{\und{*}}}_{\scriptscriptstyle{#1}}}

\newcommand{\E}[1]{\A{#1}} % edge relation
\newcommand{\Er}[1]{\Ar{#1}} %{\E^*} % edge relation - reflex-trans closure
\newcommand{\invE}[1]{\invA{#1}}
\newcommand{\invEr}[1]{\invrA{#1}} %{\invr\invE}

\newcommand{\graph}[1]{#1} %{\mbox{\sf #1}}

\newcommand{\gr}{G}

\newcommand{\hr}{{\graph{H}}}

\newcommand{\Pset}[1]{{\mathcal P}(#1)}

\newcommand{\<}{\langle}
\renewcommand{\>}{\rangle}
\newcommand{\into}{\mathrel{\rightarrow}}

\newcommand{\impl}{\into}

\renewcommand{\iff}{\leftrightarrow}
\newcommand{\Into}{\mathrel{\Rightarrow}}
\newcommand{\Impl}{\Into}

\newcommand{\Iff}{\Leftrightarrow}
\newcommand{\Ifff}{\Longleftrightarrow}

\newcommand{\by}[2]{\stackrel{#2}{#1}}

\newcommand{\eset}{\emptyset} %{\varnothing} NDJ.....
\newcommand{\false}{{\bf 0}}
\newcommand{\true}{{\bf 1}}
\newcommand{\pard}{\bot}

\newcommand{\band}{\bigwedge}
\newcommand{\bor}{\bigvee}

\newcommand{\und}[1]{\underline{#1}}

\newcommand{\nand}{\mbox{\sc nand}}
\newcommand{\ors}{\mbox{\sc or}}
\newcommand{\nandp}[2]{\ovr{#1#2}} %{\ovr{#1}\!\ \ovr{#2}}
\newcommand{\ssetminus}{\mathop{\setminus\!\!\setminus}}

\newcommand{\prr}[2]{\mbox{$\displaystyle{{#1}\over{#2}}$}}
\newcommand{\prline}[1]{\mbox{$\displaystyle{\begin{array}{c}\ \\ {#1}\end{array}}$}}
\newcommand{\pr}[2]{\vdash_{\!\!\!\tiny{#1}}^{\tiny{#2}}}
\newcommand{\pro}{\pr{}{}}

\newcommand{\noo}[1]{}

\newcommand{\Ker}[1]{Ker(#1)}

\newcommand{\sols}[1]{Ker(#1)}

\newcommand{\tr}[1]{#1^\true}
\newcommand{\fa}[1]{#1^\false}

\newcommand{\thr}{\mbox{$\mathcal T$}}
\newcommand{\cthr}{\mbox{$\mathcal C$}}
\newcommand{\RES}{\mbox{RES}}

\newcommand\subsetsim{\mathrel{%
  \ooalign{\raise0.2ex\hbox{$\subset$}\cr\hidewidth\raise-0.8ex\hbox{\scalebox{0.8}{$\sim$}}\hidewidth\cr}}}
\newcommand\supsetsim{\mathrel{%
  \ooalign{\raise0.2ex\hbox{$\supset$}\cr\hidewidth\raise-0.8ex\hbox{\scalebox{0.8}{$\sim$}}\hidewidth\cr}}}

\newcommand{\ourGr}{{\mathbf{D}}}

\title{Paraconsistency, resolution and relevance}
\author{Micha{\l} Walicki and Sjur Dyrkolbotn}
\begin{document} 
\maketitle
%\begin{frontmatter}
%runtitle{Paraconsistency, resolution and relevance}
%\author{\fnms{Micha{\l}} \snm{Walicki} \ead[label=e1]{michal@uib.no}}
%\address{Department of Informatics\\
%University of Bergen, Norway\\ \printead{e1}}
%%
%\and\ \ 
%\author{\fnms{Sjur} \snm{Dyrkolbotn} \ead[label=e2]{sdy@hvl.no}}
%\address{Department of Civil Engineering \\ Western Norway University of Applied Sciences\\ \printead{e2}}

%\runauthor{M.~Walicki and S.~Dyrkolbotn}

\hyphenation{con-di-tion-ing be-com-ing be-comes ac-cord-ing}

\begin{abstract}
  Digraphs provide an alternative syntax for propositional logic, with
  digraph kernels corresponding to classical models.
  Semikernels generalize kernels and we identify a subset of well-behaved semikernels that 
  provides nontrivial models for inconsistent theories,
  specializing to the classical semantics for the consistent ones.
  Direct (instead of refutational) reasoning with classical resolution
  is sound and complete for this semantics, when augmented with a
  specific weakening which, in particular, excludes \emph{Ex Falso}.
  Dropping all forms of weakening yields reasoning which also avoids typical fallacies of relevance.
\end{abstract}

%\begin{keyword}[class=AMS]
%\kwd[Primary ]{03B53} \kwd{03B47} \kwd{68T15}
%\end{keyword}
% 
%\end{frontmatter}

\section{Introduction} Numerous approaches to paraconsistency seem to
agree on one thing: modifications of the classical logic, made to
avoid explosion in the face of inconsistency, should be as limited as
possible.\noo{In the following, we provide a paraconsistent semantics
  for a system of resolution that satisfies this objective in an
  unusually strong sense:} We provide a paraconsistent semantics and
reasoning satisfying this objective in an unusually strong sense:
models and consequences of consistent theories are exactly their
classical models and consequences, while reasoning applies only
classical resolution. Following \cite{PatPar,ExpPow}, we give an
equivalent formulation of propositional syntax as digraphs and of
classical semantics as digraph kernels, which are generalized to
semikernels. Semikernels underlie a uniform, general concept of a
model, which gives classical models for consistent theories as a
special case. For each theory, this general concept yields a unique
set of atoms involved into inconsistency, which is empty when the
theory is consistent.

The new semantics is the main contribution of the paper.\noo{ Its
  significance is demonstrated by our proof that it has a sound and
  complete reasoning system, which is simply classical resolution.}
Its significance is supported further by some informal
  justification as well as the fact that classical resolution provides sound and
  complete reasoning.\noo{
Besides its significance is supported by the fact that 
classical resolution provides sound and complete reasoning.}
Paraconsistency of direct (instead of refutational) resolution was
applied in \cite{RIP} to infinitary logic. However, that work lacked
the semantic counterpart which is now provided. Direct resolution,
applied here to finitary (usual) propositional logic, deviates from
refutational resolution primarily by the exclusion of \emph{Ex
  Falso}. The graph syntax we use, expressed in the language of
clauses, makes \emph{Ex Falso} a special case of weakening. An
appropriate adjustment of weakening prevents then explosion from a
contradiction, allowing for its unrestricted applicability when the
theory is consistent.

Section \ref{bac} presents the background from \cite{ExpPow,RIP}, explaining
the applicability of digraphs as propositional syntax, their kernels as
classical models, and semikernels as a generalization of
kernels. Section \ref{new} presents the main contribution: a semantics 
defined in terms of well-behaved semikernels, assigning
a nonempty set of models to every theory and specializing to the
classical semantics for consistent theories. The main theorem
\ref{lemma:main} shows that every inconsistent theory has a unique set
of bad atoms, contributing to inconsistency. The consequence relation
also specializes to the classical consequence for consistent
theories. Unlike most formalisms, but in agreement with the natural
tendency of informal discourse, it disregards inconsistent parts of
statements whenever it is possible to extract from them also meaningful
elements, to which truth-values can be consistently assigned. Section \ref{logic}
shows soundness and completeness of resolution with appropriate
weakening rules. Section \ref{rel} identifies elements of relevance reasoning and their semantics, 
 arising when resolution is used without any form of weakening.

\section{Digraphs as propositional syntax}\label{bac}
A propositional formula is in graph normal form, GNF, when it has the form 
\eq{\label{GNF} x\iff \band_{i\in [n_x]}\!\!\neg y_i,}
where all $x,y_i$ are atoms (propositional variables), $n_x\in\omega$
and $[n]=\{1,...,n\}$. When $n_x=0$, the corresponding formula is
$x$. A theory is in GNF if all formulae are in GNF, and every atom of
the theory occurs exactly once unnegated, i.e., on the left of
$\iff$.\footnote{The formula $a\iff \neg b$ is in GNF but the theory
$\{a\iff \neg b\}$ is not, due to the loose $b$. Such cases can be
treated as abbreviations, here, with a fresh atom $b'$ and two
additional formulae $b\iff \neg b'$ and $b' \iff \neg b$.}
As an example we will use the formalization $\Delta$, to the
right, of the discourse to the left. The
 statement (a) requires introduction of a fresh atom
$a'$, to conform to GNF: \eq{
\begin{array}{l@{\hspace*{1.5em}}l@{\ \iff\ }l} & a' & \neg a \\
\mbox{(a) This statemenet is not false.} & a & \neg a' \\ \mbox{(b)
The previous and the next statement are false.} & b & \neg a \land
\neg c \\ \mbox{(c) The next statement is false.} & c & \neg d \\
\mbox{(d) The next statement is false.} & d & \neg e \\ \mbox{(e)
Statement (c) is false.} & e & \neg c
\end{array}
\label{EqD}
}
GNF is indeed a normal form: every theory in (infinitary)
propositional logic has an equisatisfiable one in GNF, \cite{ExpPow}
(new variables are typically needed to obtain GNF, as $a'$ above). The
classical semantics is defined in the usual way. 

GNF allows a natural reading of its equivalences as
propositional instances of the T-schema, expressing that the atom $x$
is true if and only if what it says, $\band_{i\in [n_x]}\neg y_i$, is
true. Taken in this light, a theory in GNF represents a collection of
T-schemata for the actual statements with possible, also indirect,
self-references. We therefore call a theory in GNF a \emph{discourse}
and define \emph{paradox} as an inconsistent discourse. Plausibility
of this definition, implicit in \cite{PatPar}, was argued and
exemplified in \cite{PDL,RIP} and is witnessed by the increasing
popularity of the corresponding graph representation in the analysis
of paradoxes,
\cite{Schindler2,PatPar,PDL,GrForPar,dangerousGraphs,RIP}.\footnote{GNF
finds also another application in argumentation theory in its
AI-variant following \cite{Dung}. In that context, our notion of the
(maximal) consistent subdiscourse amounts to a new semantics based on
admissible sets, whereby the acceptable extensions are the stable sets
of the maximally consistent subdiscourse of the original digraph
(argumentation framework).}

Theories in GNF and graphs are namely easily transformed into each
other.  A graph (meaning here ``directed graph'', unless qualified
otherwise) is a pair $\gr=\<G,\A\gr\>$ with $\E\gr\subseteq G\times
G$. (Overloading the notation $G$, for a graph and its set of
vertices, should not cause any confusion.)  We denote $\E\gr(x) =
\{y\in G\mid \E\gr(x,y)\}$, $\invE\gr(x) = \{y\in G\mid x\in
\E\gr(y)\}$, and extend pointwise such notation to sets, i.e.,
$\invE\gr(X) = \bigcup_{x\in X}\invE\gr(x)$, etc. $\Er\gr/\invEr\gr$
denote reflexive, transitive closure of $\E\gr/\invE\gr$.

A GNF theory $\Gamma$ gives a graph $\gr$ with all atoms as vertices
and edges from every $x$ on the left side of its GNF formula to each
$y_i$ on its righ side, i.e., $\E{\gr}= \{\<x,y_i\>\mid x\in G, i\in
[n_x]\}$. The graph for $\Delta$ from \refp{EqD} is: \eq{
\ourGr:\draw{a' \ar@<2pt>[r] & a\ar@<2pt>[l] & b\ar[l]\ar[r] &
c\ar[r]& d\ar[r] & e.\ar@/_/[ll]}
\label{graph}
} 
Conversely, the theory of a graph $\<G,\E{}\>$ is \(\thr(\gr) = \{x
\iff \band_{y\in \E{}(x)}\hspace*{-.1em} \neg y\mid x\in G\}.\) (When
$x$ is a \emph{sink}, $\E{}(x)=\eset$, this becomes $x\iff \top$,
i.e., $x$ is included in $\thr(\gr)$.) The two are inverses, so we
ignore usually the distinction between theories (in GNF) and graphs,
viewing them as alternative presentations.

The equivalence of graphs and GNF theories is not only a syntactic
transformation.  The classical models of GNF theories can be defined
equivalently as kernels of the corresponding graphs,
\cite{PatPar,ExpPow}. A \emph{kernel} of a graph $\gr$ is a subset
$K\subseteq G$ which is independent (no edges between vertices in $K$)
and absorbing its complement (every $y\in G\setminus K$ has an edge to
some $x\in K$), i.e., such that $\invE\gr(K) = G\setminus
K$. $\Ker\gr$ denotes kernels of $\gr$.

Kernel of a graph $\gr$ can be defined equivalently as a 
 2-partition $\alpha=\<\alpha^\true,\alpha^\false\>$ of the vertices $G$, 
 such that $\forall x\in G:$
 \eq{\label{sem}
\begin{array}{rrcl} 
(a)& x\in \alpha^\true &\Iff&
\forall y\in \E\gr(x): y \in \alpha^\false
\\
(b)&   x\in \alpha^\false &\Iff& \exists y\in \E\gr(x): y\in \alpha^\true.
\end{array}
}A 2-partition $\alpha$ satisfies (\ref{sem}) iff $\alpha^\true\in
 \Ker\gr$. On the other hand, satisfaction of (\ref{sem}) at every
 $x\in G$ is equivalent to the satisfaction of the respective GNF
 theory $\Gamma=\thr(\gr)$. So, for corresponding graph $\gr$ and
 theory $\Gamma$, we identify also kernels of the former and models of
 the latter. In short, graphs provide syntax for propositional logic,
 while their kernels are its classical semantics.

 $CMod(\Gamma)$ denotes classical models of $\Gamma$, each represented
 as a partition  $\alpha=\<\alpha^\true,\alpha^\false\>$ of $\gr$, 
 where $\alpha^\true/\alpha^\false$
 are atoms assigned $\true/\false$.  The classical satisfaction,
 $\models_c$, is obtained by the standard extension to complex
 formulae of the basis for atoms
 $a\in\gr: \<\alpha^\true,\alpha^\false\>\models_c a$ iff
 $a\in \alpha^\true$ and
 $\<\alpha^\true,\alpha^\false\>\models_c \neg a$ iff
 $a\in \alpha^\false$. $C$ is a classical consequence of $\Gamma$,
 $\Gamma\models_c C$ if
 $\forall \alpha\in CMod(\Gamma):\alpha\models_c C$.

 The exact correspondence between kernels of $\gr$ and models of the
 respective theory $\Gamma$ is as follows:~\footnote{Sufficient
 conditions for absence of paradox, expressed in terms of the
 properties of the graph representing the discourse, can be thus
 imported from kernel theory, as illustrated in \cite{PDL}. They
 confirm that the liar, as a minimal odd cycle, is the paradigmatic
 pattern of a finitary paradox: a finitely branching graph without odd
 cycles has a kernel. For the infinitary case, it is natural to
 conjecture that one also has to exclude some form of a Yablo
 pattern. Such a generalization is proposed in \cite{FinEnds}
 and, in an equivalent formulation, in \cite{Schindler2}. The proof of
 its special case in \cite{FinEnds} demonstrates the difficulty of the
 problem.}
 \eq{\label{eq:tight} CMod(\Gamma) = \{\<\alpha^\true,\alpha^\false\>
   \in \Pset\gr \times \Pset\gr \mid \alpha^\true \in \Ker \gr,
   \alpha^\false =\invE\gr(\alpha^\true)\}.  }
%do not distinguish sharply between them. 
%
Conditions (a) and (b) of \refp{sem} are equivalent for total
 $\alpha$ (with $\alpha^\false = G\setminus \alpha^\true$), but we
 will also consider more general structures, arising from the notion
 of a {\it semikernel}, \cite{NL71}, namely, a subset $S \subseteq G$
 satisfying:
 \eq{\label{semi}\E\gr(S) \mathrel{\accentset{(a)}\subseteq}
   \invE\gr(S)\mathrel{\accentset{(b)}\subseteq} G\setminus S.}
By (a), each $x\in\gr$ with an edge from $S$ has an edge back
 to $S$ and, by (b), $S$ is independent. $SK(\gr)$ denotes all
 semikernels of $\gr$. A semikernel $S$ is a kernel of the induced
 subgraph $\invE\gr[S] = \invE\gr(S) \cup S$. (An {\it induced}
 subgraph, or a subgraph {\it induced by} $\hr\subseteq\gr$ is $\hr =
 \<\hr, \E{\scriptscriptstyle{\gr}}\cap (\hr\times\hr)\>$.)
\begin{Example}\label{ex:main}
  The graph $\ourGr$ from \refp{graph} possess no kernel, as can be
  seen trying to assign values at $\{c,d,e\}$ conforming to
  \refp{sem}.  Its induced subgraph $\{c,d,e\}$ does not even possess
  a semikernel, but the whole graph $\ourGr$ possesses two, namely,
  $\alpha^\true=\{a\}$ and $ \beta^\true=\{a'\}$.
\end{Example}
 A semikernel can be defined
 equivalently as a 3-partition $\alpha=\<\alpha^\true, \alpha^\false, \alpha^\pard\>$ of $G$ such that
 $\forall x\in G:$
\eq{\label{semisem}\begin{array}{rrcl@{\hspace*{2em}}l} 
(a)& x\in \alpha^\true &\Rightarrow&\forall y\in \E\gr(x): y \in \alpha^\false
    & \E\gr(\alpha^\true)\subseteq \alpha^\false\\
(b)& x\in \alpha^\false &\Iff& \exists y\in \E\gr(x): y\in \alpha^\true
    & \alpha^\false = \invE\gr(\alpha^\true)\\
(c) & x\in\alpha^\pard&\Iff&x\in \gr\setminus(\alpha^\true\cup\alpha^\false) & \alpha^\pard = \gr \setminus (\alpha^\true\cup\alpha^\false).
\end{array}
} 
A 3-partition $\alpha$ satisfies \refp{semisem} iff
$\alpha^\true$ is a semikernel. $\alpha^\true$ is a kernel iff
$\alpha^\true$ is a semikernel and $\alpha^\pard = \eset$.

\section{Semantics of inconsistency}\label{new}
Models of an arbitrary theory, i.e., of a graph possibly
without any kernel, are required to satisfy three conditions which we
now briefly motivate.

In the graph $\ourGr$ from \refp{graph}, %Example \ref{ex:main}, 
the subgraph induced by $\{c,d,e\}$ has no semikernel, 
%does not admit any assignment satisfying even \refp{semisem}, 
but the subgraph induced by
$\{a',a,b\}$ -- the meaningful subdiscourse -- has two kernels: $\alpha^\true = \{a\}$ and
$\beta^\true = \{a',b\}$. The latter does not seem adequate as a model, because
it should function in the context of the whole original
theory, and not only after removal of its inconsistent part. In the
context of the whole $\ourGr$, $b$ negates not only $a$ but also
$c$, so to conform to \refp{sem}, or even just \refp{semisem},
$b\in \beta^\true$ would require $c\in\beta^\false$. %, while $c\in \beta^\pard$. 
Choosing $\alpha^\true$ instead, $b\in \alpha^\false$ complies with
\refp{semisem} since $a\in \E{\ourGr}(b) \cap \alpha^\true$.

This suggests semikernels as a semantic basis in the presence of
inconsistency: $\alpha^\true\in SK(\ourGr)$, while
$\beta^\true \not\in SK(\ourGr)$. A semikernel $\tr\alpha$ makes all
$x\in \alpha^\true$ fully justfied, in the sense that
$\E\gr(x)\subseteq \alpha^\false$. This excludes $\beta^\true$ from
possible models. For $x\in \alpha^\false$, on the other hand, it
suffices that $\E\gr(x)\cap \alpha^\true\not=\eset$, allowing other
out-neighbours of $x$ to be arbitrary -- possibly paradoxical. Such
paradoxical elements form the third part $\alpha^\pard$ of the model.

There are, however, too many semikernels. In the graph $\ourGr$, each
among $\{a\}$, $\{a'\}$ and $\eset$ (giving
$\alpha^\true=\alpha^\false=\eset$ and $\alpha^\pard=\ourGr$) is a
semikernel. Such a semantics is too liberal and we have to choose
semikernels more carefully.

To explain next restriction, it will be helpful
to consider two simple examples. ``This statement is
false and the sun is not a star'' is represented by
$F_1$:\hspace*{1.5em}\draw{f\ar@(dl,ul)\ar[r] & s.} Here, $f$
seems false, negating the true statement $s$. Indeed, $F_1$ has a
unique kernel, $\alpha^\true=\{s\}$, which yields $f\in \fa\alpha = F_1\setminus \tr\alpha$. Now,
consider ``This statement is false and %it is false that the sun is not
the sun is a star'' -- represented by
$F_2$:\hspace*{1.5em}\draw{f\ar@(dl,ul)\ar[r] & y\ar[r] & s.} Now
$f$ appears to be paradoxical, since $F_2$ is a contingent liar, ceasing to be 
paradoxical only if the sun is not a star, which it is.  The only
semikernel of $F_2$, $\alpha^\true = \{s\}$, gives
$\alpha^\false=\{y\}$, but this leads to the irresolvability of the
paradox ``at'' $f$. In short, $f$ is paradoxical because $s$
happens to be true. Hence, also $s$ is involved in the paradox at $f$,
not as a standalone atom but as a member of the formula for $y$, which
in turn features in the formula of the contingently paradoxical $f$.

The paradox ``at'' $f$ -- ``referring to'' $s$ by denying $y$ --
involves $s$ and $y$ as much as it involves self-reference. If we are
not prepared to admit this, we should hardly regard $F_1$ as
nonparadoxical, since the problematic self-reference at $f$ is
exactly the same in both $F_1$ and $F_2$.  If the truth of $s$
prevents paradox in $F_1$ then, in the same way, it contributes to it
in $F_2$.

This is not to suggest that ``the sun is a star'' is paradoxical on
its own, only that its token contributes to the paradoxical whole when
combined with the contingent liar as in $F_2$. Consistency and paradox
are genuinely holistic.  Or put differently: the token of ``the sun is
a star'' is unproblematic in $F_1$, but its token in $F_2$ {\it
  becomes} paradoxical by contributing to the appearance of the
paradox: if there were no $s$, there would be no paradox.  When trying
to repair this paradox, removing the loop at $f$ is as good as
removing $s$.

The absence of any single culprit among $\{f,y,s\}$ is just as it was
with $\{c,d,e\}$ in \refp{EqD}.  A nonobvious
informal lesson could be: if an inconsistency, occurring ``at'' some
$f$, depends on some $s$ (in the sense of $s\in \Er\gr(f)$), then $s$
is ``a part of'' this inconsistency.  Consequently, we should not rest
satisfied with an arbitrary semikernel, like $\{s\}\in SK({F_2})$.
Exactly the semikernel we choose (combined with other factors, like
the loop at $f$) can be the reason for the inconsistency, which could
be possibly prevented by another choice. A satisfactory semikernel $S$
should not contribute to any inconsistency occuring above it in the
graph. Put precisely, a model is not only a semikernel but an
$\invE\gr$-closed semikernel, i.e.:
\eq{S\in SK(\gr): \invE\gr(\invE\gr[S])\subseteq \invE\gr[S].\label{M2}}
This views all $\{f,y,s\}$ as inconsistent (contributing to the inconsistency) in $F_2$,
making the empty semikernel %, $\alpha^\true=\eset$, 
the only
interpretation. For $\ourGr$ from \refp{graph}, 
 the semikernel $\{a'\}$ is rejected, as it is not
$\invE\gr$-closed, leaving only  $\{a\}$. 

The above condition still admits the empty semikernel, even when there
are nonempty ones. To avoid this we require $\invE\gr[S]$ to be
maximal: %, i.e., such that:
\eq{\label{M3}
%\forall R\subseteq \gr: R \rm{\ satsfies\ }\refp{M2} \Impl \invE\gr[R]\subseteq\invE\gr[S].
\forall R\subseteq \gr: R \rm{\ satsfies\ }\refp{M2} \Impl \invE\gr[S]\not\subset\invE\gr[R].
}
%
%Note that the negated inclusion is strict. 
This condition can be seen
as a minimization of inconsistencies, typical for preferential models,
like LP$^{\rm m}$ and many other examples.  But here most of such minimization
is done by the two earlier conditions; this one excludes only specific
degenerate cases.  We thus obtain the main definition: models of a
graph $\gr$ (or its theory $\Gamma$) are $\invE\gr$-closed, maxi\-mal
semikernels:
\eq{
\begin{array}{r@{\hspace*{-.01cm}}l}
Mod(\gr)  = & \{ S\in SK(\gr)\mid 
\invE\gr(\invE\gr[S])\subseteq \invE\gr[S]\  \land\  \forall R\subseteq\gr: 
\\[.5ex]
& \  
R\in SK(\gr)\! \land 
\invE\gr(\invE\gr[R])\subseteq\! \invE\gr[R] \Impl 
\invE\gr[S]\not\subset\invE\gr[R] \}
\end{array} 
\label{Mod}}
The consequence relation $\models$ generalizes $\models_c$ to
3-partitions $\alpha_S=\<\tr\alpha_S, \fa\alpha_S, \alpha^\pard_S\>$
$= \<S, \invE\gr(S), \gr\setminus \invE\gr[S]\>$, arising from all 
independent $S\subseteq\gr$, in particular, models of $\gr:
\Gamma\models C$ if $\forall S\in Mod(G):\alpha_S\models C$.  As
formulae, we use clauses and let $A\ovr B$ denote a clause with atoms
$A$ and negated atoms $B$.  Such a clause is satisfied by a
3-partition $\alpha$ according to the following rule:
%The satisfaction relation is defined by
%
\eq{\label{parasem} \<\alpha^\true,\alpha^\false,\alpha^\pard\>\models
  A\ovr B\ \Iff\ A\cap \alpha^\true \not=\eset \lor B\cap \alpha^\false
  \not=\eset \lor (AB\subseteq \alpha^\pard\not=\eset). } 
The third disjunct may appear less intuitive than the first two, but
we comment it below. When $\alpha^\pard=\eset$, this reduces to the
classical satisfaction and consequence with respect to $CMod(\gr)$, 
%as defined before 
cf. \refp{eq:tight}.  In the extreme case of a $\Gamma$
where all atoms are involved in inconsistency, there is only one model
$\<\eset,\eset,\gr\>$, arising from the empty semikernel, and
satisfying every clause (over the atoms $\gr$).  On the other hand,
the empty theory has the empty graph, with only one, empty kernel,
giving the only model $\<\eset,\eset,\eset\>$, which does not satisfy
any formula.  This deviation from classical logic is only a question
of preference, since for this special case we could define the models
in the classical way. However, since our logic is a logic of
consequences rather than of tautologies, it appears plausible that
nothing follows from saying nothing.

According to \refp{parasem}, an $\alpha$ satisfies a clause $C$ either
for some ``healthy'' reasons, when some of its literals are true, or
because the clause is ``completely nonsensical'', with all atoms
involved into inconsistency. In natural discourse, we tend to focus on
its meaningful parts, simply ignoring occasional nonsense. A statement
``The sun is a planet; or else this (part of this) statement is
false.'' may be judged nonsensical (as it would be in SK or LP, where
$\false\lor \pard =\pard$). But if we grant the interlocutor the
benefit of doubt and are willing to ignore the partial nonsense, we
can also say that it is false, since so is its meaningful part. When,
however, unable to discern any sense whatsoever, like in the liar or
in (c)-(d)-(e) from \refp{EqD}, we ``accept'' the claim as much as we
``accept'' its negation.  Relation \refp{parasem} can be read as such
an acceptance which treats clauses containing healthy literals
according to these literals, ignoring the nonsensical part. Faced with
a complete nonsense, however, it becomes as confused as we are when,
in the face of the liar, we find it equally (im)plausible to accept
its ``truth'' and its ``falsehood''. This does not imply any semantic
dialetheism, since nonsensical atoms are excluded from the healthy
considerations and relegated to $\alpha^\pard$. (They can be seen as
gluts, since both the atom and its negation are provable, but also as
gaps, being irrelevant for the value of clauses containing also
healthy literals.)  A significant point is that this acceptance
relation is not used for defining the models, which are chosen using
\refp{Mod}, but only for determining their consequences. As it
happens, the members of $Mod(\gr)$ do satisfy the graph's theory
according to \refp{parasem}, but they need not be all 3-partitions
doing this.

Typically, $Mod(\gr)$ contains nontrivial models also for inconsistent
$\Gamma$. In fact, these are classical models of the appropriate
subgraph.  Given an $\alpha^\true\in Mod(\gr)$, and projecting away
the third component from its partition
$\<\alpha^\true,\invE\gr(\alpha^\true),\gr\setminus
\invE\gr[\alpha^\true]\>$, leaves
$\<\alpha^\true,\invE\gr(\alpha^\true)\>\in
CMod(\invE\gr[\alpha^\true])$,
i.e., $\alpha^\true \in Ker(\invE\gr[\alpha^\true])$ -- a classical
model of the theory for the induced subgraph
$\invE\gr[\alpha^\true]$.  Interestingly, each two of such models
$\tr\alpha,\tr\beta\in Mod(\gr)$ classify the same vertices as
inconsistent, assigning boolean values to the same subset of $\gr$,
namely, $\invE\gr[\tr\alpha] = \invE\gr[\tr\beta]$.

Proving this will take the rest of this section and requires some
preliminary observations.\noo{ First, each independent $S\subseteq \gr$
determines a 3-partition $\alpha_S$ of $\gr$ given by:
$\alpha_S^\true = S$, $\alpha_S^\false = \invE\gr(S)$ and
$\alpha_S^\pard = G\setminus (\alpha^\true\cup \alpha^\false)$.} When
$S\in SK(\gr)$, $\alpha_S=\<S,\invE\gr(S),\gr\setminus \invE\gr[S]\>$ satisfies conditions (a) and (b) from
\refp{semisem}, repeated below: \vspace*{-2.5ex}

\eq{\label{mps}\begin{array}{rrcl@{\hspace*{2em}}l} 
(a)& x\in\alpha^\true &\Rightarrow& \forall y\in \E\gr(x): y \in\alpha^\false 
    & \E\gr(\alpha^\true)\subseteq \alpha^\false\\
(b)& x\in \alpha^\false &\Iff& \exists y\in \E\gr(x): y\in \alpha^\true
    & \alpha^\false = \invE\gr(\alpha^\true)\\
(c)&   x\in \alpha^\pard &\Iff& \forall y\in \E\gr(x): y\in \alpha^\pard 
    & \alpha^\pard = \E\gr[\alpha^\pard]
\end{array} 
}
Condition (c) in \refp{semisem} was a mere definition of
$\alpha^\pard = \gr \setminus (\tr\alpha\cup\fa\alpha)$, while here it
expresses that $S$ is $\invE\gr$-closed, since its complement
$\alpha^\pard$ is $\E\gr$-closed. \pagebreak[3]Obviously,
$\invE\gr(\invE\gr[S]) \subseteq (\invE\gr[S])$ implies
$\alpha_S^\pard\supseteq\E\gr[\alpha_S^\pard]$. The opposite
implication of (c) follows then from (a) and (b), while
$\alpha_S^\pard\subseteq\E\gr[\alpha_S^\pard]$ holds by definition
$\E\gr[X]=X\cup \E\gr(X)$.  (Unlike in (a) and (b), the two
formulations in (c) are not equivalent.) Hence for every
$\invE\gr$-closed $S\in SK(\gr)$, $\alpha_S$ satisfies \refp{mps}, and
$pSK(\gr)$ denotes all such 3-partitions.

Conversely, every 3-partition $\alpha\in pSK(\gr)$  satisfies 
 \refp{semisem}, so that $\alpha^\true\in SK(\gr)$. 
% Furthermore, $\alpha$ satisfies also 
Such an $\alpha$ satisfies also the following closure property:
\eq{\label{ref:closure}
\begin{array}{ll}
(i) & \E\gr(\alpha^\pard) \subseteq \alpha^\pard \\
(ii) & \invE\gr(\alpha^\true \cup \alpha^\false) \subseteq \alpha^\true \cup \alpha^\false
\end{array}
} 
Point (i) follows from condition (c) of \refp{mps} while point (ii) is
equivalent to (i), since
$\alpha^\pard = \gr \setminus (\alpha^\false\cup\alpha^\true)$. When
$\alpha^\true\in SK(\gr)$, (ii) is the condition \refp{M2}.  Thus pSK
partitions correspond exactly to $\invE\gr$-closed semikernels,
\refp{M2}, so that $Mod(\gr)$ correspond exactly to \emph{maximal} pSK
partitions, namely: %$mpSK(\gr)$, namely:
\eq{\label{mpSK} \alpha\in mpSK(\gr)  \Iff \alpha\in pSK(\gr) \land 
 \forall \beta \in pSK(\gr): 
  \alpha^\false\cup\alpha^\true \not\subset
\beta^\false\cup \beta^\true. }
In addition, the following fact is used in the proof of the next, crucial 
lemma.
\begin{Fact}\label{fact:comb}
For every graph $\gr$,
\begin{enumerate}
\item If $T\subseteq S\in SK(\gr)$  then $(S \cap \Er\gr(T)) \in SK(\gr)$.
\item If $S\in SK(\gr)$, $T\in SK(\gr)$  and
  $\invE\gr(S) \subseteq G \setminus T$, then $(S \cup T)\in SK(\gr)$. 
\end{enumerate}
\end{Fact}
\begin{Proof}
  1. For any
  $x\in \E\gr(t)\subseteq \E\gr(T) \subseteq \E\gr(S)\subseteq
  \invE\gr(S)$,
  and $s\in S$ with $x\in \invE\gr(s)$, we have
  $s\in \Er\gr(t)\subseteq\Er\gr(T)$, i.e.,
  $x\in \invE\gr(S\cap \Er\gr(T))$. This gives the first inclusion
  below, while the second one follows since $S\in SK(\gr)$:

  $\E\gr(S \cap \Er\gr(T)) \subseteq \invE\gr (S \cap \Er\gr(T))
  \subseteq \gr\setminus S \subseteq \gr \setminus (S \cap
  \Er\gr(T))$.
  \\
  2.
  $\E\gr(S\cup T) \subseteq \E\gr(S)\cup \E\gr(T) \subseteq
  \invE\gr(S)\cup \invE\gr(T) = \invE\gr(S\cup T)$.
  For the next inclusion, we note that
  $\invE\gr(S)\subseteq\gr\setminus T$ implies here also the dual
  $\invE\gr(T)\subseteq\gr\setminus S$, for if for some
  $t\in T,s\in S: t\in \E\gr(s)$, then $t\in \invE\gr(S)$ since
  $S\in SK(\gr)$.  Hence

  $\invE\gr(S)\cup\invE\gr(T) \subseteq ((\gr\setminus S)\cap
  (\gr\setminus T)) \cup ((\gr\setminus T) \cap (\gr\setminus S)) =
  \gr \setminus (S\cup T)$.
\end{Proof}
Consequently, distinct $\invE\gr$-closed semikernels, %\refp{M2}/\refp{mps},
disagreeing on at least one paradoxical element, can be combined as in
the proof of the following lemma.
\begin{Lemma}\label{lemma:main0}
For all graphs $\gr:$ 

$\forall \alpha,\beta\in pSK(\gr)\ \exists \gamma\in pSK(\gr): 
\beta^\true\cup\beta^\false \not\subseteq \alpha^\true\cup\alpha^\false \Impl 
\alpha^\true\cup\alpha^\false \subset \gamma^\true\cup\gamma^\false$. 
\end{Lemma}
\begin{Proof} For arbitrary $\alpha,\beta \in pSK(\gr)$ with
  $\beta^\true\cup\beta^\false \not\subseteq
  \alpha^\true\cup\alpha^\false$,
  we have that
  $(\beta^\true\cup\beta^\false) \cap \alpha^\pard \not=\eset$, so
  define $Q = \beta^\true \cap \alpha^\pard$ and
  $R = \beta^\false \cap \alpha^\pard$. We show that
  $S = \alpha^\true \cup Q$, with $Q\not=\eset$, is an
  $\invE\gr$-closed semikernel, i.e., the desired $\gamma^\true$.
\begin{itemize}
\item [(a)] $R \subseteq \E\gr(Q)$, by $\beta^\true\in SK(\gr)$ and
  \ref{ref:closure}.(i) -- hence also $Q\not=\eset$. 
\item [(b)] $\Er\gr(Q) \subseteq \Er\gr(\alpha^\pard)\subseteq \alpha^\pard$, 
by \ref{ref:closure}.(i).
\item [(c)] $Q$ is a semikernel, because 
 $\beta^\true\cap \Er\gr(Q)\in SK(\gr)$ by Fact \ref{fact:comb}.(1), 
while  $Q=\beta^\true\cap \Er\gr(Q)$ 
by (b): 
$\beta^\true \cap \alpha^\pard \subseteq \beta^\true\cap \Er\gr(Q)\subseteq \beta^\true\cap\alpha^\pard$.
\item [(d)] $\invE\gr(Q) \subseteq G \setminus \alpha^\true$, by
  $Q \subseteq \alpha^\pard$ and $\alpha^\true\in SK(\gr)$ (so that
  $\invE\gr(\alpha^\pard) \cap \alpha^\true = \eset$).
\item [(e)] $S\in SK(\gr)$, by Fact \ref{fact:comb}.(2) (applicable by (c)-(d) above).
\item [(f)] $S$ is $\invE\gr$-closed, i.e., $\invE\gr(\invE\gr[S]) \subseteq \invE\gr[S]$. If
  $x \in S \subseteq \invE\gr[S]$, then trivially
  $\invE\gr(x) \subseteq \invE\gr[S]$. If $x \in \invE\gr(S)$, we have two
  cases. \begin{itemize}
  \item [(i)] $x \in \invE\gr(\alpha^\true)$. 
Since $\alpha$ is a pSK partition, \refp{mps}: 
    $\invE\gr(x) \cap \alpha^\pard = \eset$. Since
    $G \setminus \invE\gr[S] \subseteq \alpha^\pard:
    \invE\gr(x) \subseteq \invE\gr[S]$, as desired.

\item [(ii)] $x \in \invE\gr(Q)$.  $\alpha$ satisfies \refp{mps} and
$\beta$ \refp{semisem}, so $\invE\gr(Q) \subseteq \beta^\false \cap
(\alpha^\pard \cup \alpha^\false)$.

  If $x\in \alpha^\false$ then $\invE\gr(x)\subseteq
  \alpha^\true\cup\alpha^\false$, because $\alpha$ satisfies
  \refp{mps}. Hence $\invE\gr(x)\subseteq \invE\gr[S]$.  If $x\in
  \beta^\false\cap \alpha^\pard$ then $\invE\gr(x)\subseteq
  (\beta^\true\cup\beta^\false)\cap
  (\alpha^\pard\cup\alpha^\false)$. For any $y\in \invE\gr(x)$:

  $y \in (\beta^\true\cup\beta^\false) \cap\alpha^\false \subseteq
  \invE\gr(\alpha^\true)\subseteq \invE\gr(S)$, or

  $y \in \beta^\true \cap\alpha^\pard = Q \subseteq \invE\gr[S]$, or 

  $y \in \beta^\false \cap\alpha^\pard$ -- then
$\exists z\in \beta^\true: y\in \invE\gr(z)$, since $\beta^\true\in SK(\gr)$.
 Since $z\in \E\gr(y)$, so $z\in \alpha^\pard$ by
\ref{ref:closure}.(i), which means that $z\in Q$ so that
$y\in \invE\gr(Q)\subseteq \invE\gr(S)$.
\end{itemize}\vspace*{-3.5ex}
\end{itemize}%\vspace*{-2ex}
\end{Proof}
So, if $\alpha,\beta\in mpSK(\gr)$ while
$\alpha^\pard\not=\beta^\pard$, then
$\beta^\true\cup\beta^\false \not=\alpha^\true\cup\alpha^\false$, 
in particular,
$\beta^\true\cup\beta^\false
\not\subseteq\alpha^\true\cup\alpha^\false$
(since $\subset$ would contradict \refp{mpSK} for $\beta$). %\in mpSK(\gr)$). 
%By the above, t
There is
then
$\gamma\in pSK(\gr): \alpha^\true\cup\alpha^\false \subset
\gamma^\true\cup\gamma^\false$,
contradicting \refp{mpSK} for $\alpha$. We thus obtain:
%, the domain of the
%boolean assignments in all $Mod(\gr)$ is determined uniquely:
%
\begin{Theorem}\label{lemma:main}
For all graphs $\gr$ and all $S,R\in Mod(\gr): \invE\gr[S]=\invE\gr[R]$. 
\end{Theorem}
Since every $S\in SK(\gr)$ 
%Every semikernel $S$ 
is a kernel of the subgraph induced by $\invE\gr[S]$,\noo{but our
semantics \refp{Mod} satisfies a much stronger property: all
$Mod(\gr)$ are semikernels of $\gr$ which are also kernels of the same
induced subgraph of $\gr$, namely, $\invE\gr[S]$, for any maximal
$\invE\gr$-closed semikernel $S$.} our $Mod(\gr)$ satisfy quite a
strog property: they are namely kernels of one specific subgraph of
$\gr$, given by $\invE\gr[S]$, for any maximal $\invE\gr$-closed
semikernel $S$. When the theory is consistent, $\invE\gr[S]=\gr$ and 
%For the special case of a consistent theory, 
$Mod(\gr) = \Ker\gr$, i.e., the models 
are exactly the classical ones.

The induced subgraph $\invE\gr[S]$, for any $S\in Mod(\gr)$, gives the
\emph{maximal consistent subdiscourse} of $\gr$.  Since the graph
provides the syntax of a theory, a (typically induced) subgraph
corresponds to a kind of subtheory, referred to as a
subdiscourse. This concept differs from a subtheory seen as a subset
of the formulae. In the graph $\ourGr$ from \refp{graph}, $Mod(\ourGr)
\subseteq \Ker\hr$, where $\hr$ is the induced subgraph
\draw{a'\ar@<2pt>[r]& a\ar@<2pt>[l] & b\ar[l],} with the theory
$\thr(\hr)= \{b\iff \neg a$, $a\iff \neg a', a'\iff \neg a\}$. Its
formula $b\iff\neg a$ does not occur in the original theory $\Delta =
\thr(\ourGr)$ from \refp{EqD}, which has instead $b\iff \neg a \land
\neg c$.  A subdiscourse, as an induced subgraph, amounts not only to
a subset of the formulae but also, for each retained formula, possibly
only a subset of the (negated) atoms under the conjunction in its
right side.

Definition \refp{Mod} chooses as $Mod(\ourGr)$ only $\{a\}\in
\Ker\hr$, making $a=\true$ and $b=\false=f$.  The other kernel
$\{a',b\}\in\Ker\hr$ is not a semikernel of $\ourGr$, while the other
semikernel $\{a'\}$ of $\ourGr$ is not $\invE{\ourGr}$-closed. 
%This effect 
The exact subset of kernels of $\hr$
constituting the models of $\ourGr$ can be captured as the
classical models $CMod$ from \refp{eq:tight}
-- not, however, of $\hr$ but of its
appropriate modification, namely, as the kernels of
\draw{a'\ar@<2pt>[r]& a\ar@<2pt>[l] &
b\ar[l]\ar@(dr,ur)&\ .} The new loop at $b$ keeps track of
the edge $b\into c$, which disappeared in $\hr$ but prevents $b$ from
being $\true$. After all, in the original $\ourGr$, $b$ negates $c$
which is not unproblematically $\false$.\noo{  As shown in Theorem \ref{not-main}, 
this is a general phenomenon. Models of
every discourse are the classical models of such a modification of its 
maximal consistent subdiscourse. }
Theorem \ref{not-main} below shows that models of every discourse are
exactly the classical models of such a modification of its maximal consistent
subdiscourse.

\section{Reasoning}\label{logic}
The system \RES\ consists of the axiom $\pro a\ovr
a$, for every atom $a$, and the resolution rule,
\prr{\Gamma\pro Aa \ \ \ \ \Gamma\pro B\ovr a}{\Gamma\pro AB}. Clauses
are obtained from the two implications of GNF formulae
(\ref{GNF}). For each $x\in G$, they are of two kinds: %\vspace*{.5ex}

\begin{tabular}{rl} 
\ors-clause:&  $x\lor \bor_{i\in [n_x]} y_i$, written as $xy_1y_2...y_n$
\\[.5ex]
\nand-clauses:& $\neg x \lor \neg y_i$, for every $i\in [n_x]$, denoted $\nandp {x}{y_i}$. 
\end{tabular}
\\[1ex]
In terms of a graph $\gr$, its clausal theory $\cthr(\gr)$
contains, for every $x\in G$, the \ors-clause $\E\gr[x] = \{x\}\cup
\E\gr(x)$ and for every $y\in \E\gr(x)$, the \nand-clause $\nandp
xy$. For the graph $\ourGr$ from \refp{graph}, its
clausal theory is: 

$\cthr(\ourGr) = \{aa', \nandp aa', bac, \nandp ba, \nandp bc, cd,
\nandp cd, de, \nandp de, ec, \nandp ec\}$.  
\\
We treat clauses as (finite) sets of literals, with overbars marking
 the negative ones, i.e., we write $A\ovr B$ for a clause with atoms
 $A$ and negated atoms $B$. Initial uppercase letters $A,B,C...$ denote typically
 arbitrary clauses, which is sometimes marked by $C\subseteq \any
 \gr$, where $\any\gr = \gr \cup \{\ovr x\mid x\in \gr\}$.  For a
 clause $C$, $C^-$ denotes the set of unary clauses with its
 complementary literals.

Of primary interest to us are graphs (GNF theories) but several
results hold for arbitrary clausal theories (sets of finite clauses). By ``every
$\Gamma$'' we mean such theories.  (For graphs, the axiom schema is not needed, being provable for every vertex with outgoing edges, e.g., in $\cthr(\ourGr)$, resolving $bac$ with $\nandp ba$ and $\nandp bc$ yields $b\ovr b$, etc.) 
The following gathers some relevant
facts about resolution:
\begin{Fact}\label{RES}For every $\Gamma$ over atoms $\gr$ and a clause $C\subseteq\any\gr :$
\begin{enumerate}%\MyLPar
\item\label{snd}
$\Gamma\pro {C} \Impl \Gamma\models_c {C}$,
\item\label{ref}
$CMod(\Gamma) = \eset \Iff \Gamma \pro\ecl$,
\item\label{refCompl}
$\Gamma\models_c C \Iff \Gamma,{C}^-\pro \ecl$,
\item\label{wCompl}\label{cor:compl+}
$\Gamma\models_c {C} \Iff \exists {B}\subseteq {C}: \Gamma\pro {B}$,
\item\label{xnotx}
$\Gamma \pro \ecl \Iff \exists a\in \gr: \Gamma\pro a \land \Gamma\pro \ovr a$\hfill denoted $\Gamma\pro \pard(a)$,
\item\label{aux}Denoting $\RES(\Gamma) = \{C\mid \Gamma\pro C\}$: 

$\RES(\Gamma,A^-) = \RES(\Gamma)\cup A^- \cup \{P\setminus B\mid \Gamma\pro P \ {\rm and}\  B\subseteq A\}.$
\end{enumerate}
\end{Fact} 
%Point \ref{refCompl} is refutational (soundness and) completeness of
%resolution, while \ref{wCompl} its weaker and more specific
%completeness for direct (not refutational) reasoning.
%
For diagnosing inconsistency of $\Gamma$, pinpointing the problem to
specific atoms is not necessary, and it suffices that
$\Gamma\pro\ecl$, as guaranteed by point \ref{ref}. This point implies also
refutational completeness with respect to all classical consequences,
\ref{refCompl}. But we consider instead only direct (not refutational)
derivability, i.e., we ask if $\Gamma\pro  C$, instead of
$\Gamma, C^-\pro \ecl$. This gives weakened completeness, \ref{wCompl},
which could be repaired by adding the weakening rule. But its absence, 
 and the consequent inadmissibility of \emph{Ex Falso},
arise now as virtues rather than vices. They give a paraconsistent
ability to contain paradox and reason about the subdiscourse
unaffected by it.

As a simple example, for $\Gamma = \{x,\ovr x,s\}$, we have
 $\Gamma\pro\ecl$ but also $\Gamma\not\pro \ovr s$.  Its graph --
 $\draw{x\ar@(dr,ur) & \ \ \ s}$ -- justifies this: the liar $x$ is
 in no way ``connected'' to $s$.  This is the essence of the
 phenomenon, which we now describe in more detail.
\begin{Example}\label{ex:nonexplode}
The closure of $F_2:$\hspace*{1.7em}\draw{\ar@(dl,ul) f \ar[r] & y\ar[r] & s,} that is, 
of its clausal theory 
$\Gamma_2 = \{\ovr f, fy, ys, \nandp ys, s\}$ 
contains, besides $\ecl$, all literals. 
\\
The clausal theory $\cthr(\ourGr) = \{aa', \nandp aa', bac,
\nandp ba, \nandp bc, cd, \nandp cd, de, \nandp de, ec, \nandp ec\}$
 is provably paradoxical, but neither $b$, $a'$ nor
$\ovr a$ are provable. Its deductive closure contains $\pard(x)$ for
each $x\in \{c,d,e\}$ and, besides that, only $\ovr b$, $a$ and $\ovr
a'$. It determines thus the only member of $Mod(\ourGr)$.
\end{Example}
This is no coincidence -- \RES\ derives clauses satisfied by all
$Mod(\gr)$ and, extended with appropriate weakening, exactly these
clauses, but proving this will take the rest of this section.  First,
we register soundness of \RES\ for all 3-partitions.
%of the graph.
%
\begin{Fact}\label{sound}
For every 3-partition $\alpha$ of $\gr$, every $a\in \gr$ and $A,B\subseteq\any\gr:$ 

$\alpha\models a\ovr a$ and $\alpha\models Aa\land \alpha\models B\ovr a \Impl \alpha\models AB$.
%\RES\ is sound with respect to the relation $\models$ defined in \refp{parasem}. 
\end{Fact}
\begin{Proof}
Let $\alpha=\<\alpha^\true,\alpha^\false,\alpha^\pard\>$ be an
arbitrary 3-partition of $G$. Obviously, $\alpha\models a\ovr a$, for if $a\not\in \tr\alpha\cup\fa\alpha$, then $a\in \alpha^\pard$. 
Assume that $\alpha\models Aa$ and
$\alpha\models B\ovr a$.  If $AB\subseteq \alpha^\pard$, then
$\alpha\models AB$, also when $AB=\ecl$, since then both
$\alpha\models a$ and $\alpha\models \ovr a$, which imply $a\in
\alpha^\pard$.  Assuming $AB\not\subseteq \alpha^\pard$, either
$A\not\subseteq \alpha^\pard$ or $B\not\subseteq \alpha^\pard$. Wlog.,
assume $A\not\subseteq \alpha^\pard$. If $B\subseteq \alpha^\pard$,
then $a\in \fa\alpha\cup\alpha^\pard$ because $\alpha\models B\ovr
a$. Then $A$ contains a literal which witnesses to its truth (positive
in $\tr\alpha$ or negative in $\fa\alpha$), and to the truth of the
conclusion. If $B\not\subseteq \alpha^\pard$, then either $a\in
\fa\alpha$, and the conclusion follows as in the previous case, or $B$
contains a literal witnessing to its truth, and to the truth of the
conclusion.
\end{Proof} 
Consequently, % and unsurprisingly, 
our semantics %\refp{Mod}, \refp{parasem} 
agrees with the classical one as
to which theories count as inconsistent and \RES\ proves
inconsistency for exactly these theories:
%\eq{\begin{array}{rrcl}
%(a) & \Gamma\models\ecl & \Iff & \Gamma\models_c\ecl \\
%(b) & \Gamma \models \ecl & \Iff & \Gamma\pro \ecl  \end{array}}

(a) $\Gamma\models\ecl \Iff \Gamma\models_c\ecl$

(b) $\Gamma \models \ecl \Iff \Gamma\pro \ecl$.
\\
  (a) follows by $\Gamma\models_c\ecl \Iff CMod(\Gamma)=\eset
  \by\Impl{\ref{RES}.\ref{ref}}\Gamma\pro\ecl \by\Impl{\ref{sound}}
  \Gamma\models\ecl$. 
  Conversely, if $\Gamma\models\ecl$ then for every
  $\alpha\in Mod(\Gamma): \alpha^\pard \not=\eset$, i.e., 
  $CMod(\Gamma)=\eset$, so $\Gamma\models_c\ecl$.
  Combining (a) with 
  $\Gamma\pro \ecl\Iff\Gamma \models_c \ecl$ from Fact \ref{RES}.\refp{ref} yields (b). 

The rest of this section is concerned primarily with the situations
when $\Gamma\pro\ecl$ and $\exists x\in \gr:
\Gamma\not\pro \pard(x)$. The constructions and results are general,
but they trivialize when one of these conditions is violated.  Given
an arbitrary $\Gamma$, we construct $\okk\Gamma$ -- the maximal
consistent subdiscourse, with the additional requirement on its border
vertices, which refer to the inconsistent elements.  ($\okk\Gamma$ is
empty if all atoms are inconsistent and coincides with $\Gamma$, if
none is.)  The classical models of $\okk\Gamma$ turn out to be the
models of $\Gamma$, Theorem \ref{not-main}. This leads to the
completeness of \RES, where every clause satisfied by all $Mod(\gr)$
has a nonempty provable witness, Corollary \ref{cor:main}. Augmenting
\RES\ with appropriate weakening yields then a strongly complete
reasoning system.  Some technicalities below, originating from
\cite{RIP}, are adjusted to the present context and repeated to make
the paper self-contained. The main results, Theorem \ref{not-main} and
Corollary \ref{cor:main}, are new.

For a clausal theory $\Gamma\subseteq \Pset{\any \gr}$ and $X\subseteq \gr$, the
operation
\begin{center}
\(\Gamma\ssetminus X = \{C\setminus \any X\mid C\in\Gamma\}\setminus
\{\ecl\}\) 
\end{center}
removes all literals over atoms $X$ from all clauses of $\Gamma$, 
removing also the empty clause, if it appears. It satisfies the following
important property, relating consequences of $\Gamma$ to consequences
of $\Gamma \ssetminus X$.
\begin{Lemma}\label{lemma}
  For each $\Gamma$ and 
  $A \not \subseteq \any X: \Gamma \pro A \Rightarrow \exists B
  \subseteq A \setminus \any X: \Gamma \ssetminus X \pro B$. 
\end{Lemma}
\begin{Proof}
If $\Gamma\ssetminus X\pro\ecl$ the claim follows, so we assume that
this is not the case and proceed by induction on the length of the
proof $\Gamma\pro A$, with axioms introducing ${A\setminus
\any X}$ instead of $A$. A step \prr{\Gamma\pro A_1a\ \ \ \
\Gamma\pro A_2\ovr a}{\Gamma\pro A_1A_2}, where $A = A_1A_2$, has
by IH the corresponding proofs $\Gamma\ssetminus X\pro
(A_1a)\setminus \any X$ and $\Gamma\ssetminus X\pro (A_1a)\setminus \any X$. If
$a\in \any X$, either of these proofs can serve as the conclusion. 
Otherwise, conclusion follows by the deduction \prr{\Gamma\pro (A_1\setminus \any X)a\ \
\ \ \Gamma\pro (A_2\setminus \any X)\ovr a}{\Gamma\pro A_1A_2\setminus \any X}.
\end{Proof}
Hence, if $\Gamma\pro C$ and $C$ contains at least one literal not in
$\any X$, Then $\Gamma\ssetminus X\pro C'$ for some nonempty
$C'\subseteq C$: removing literals from $\Gamma$ using $\ssetminus$,
results at most in sharpening the information about the remaining
atoms.  Let us denote:

$G^\pard = \{x\in G\mid \Gamma\pro x \land \Gamma\pro \ovr x\}$

$\okk\Gamma = \Gamma\ssetminus G^\pard = \{C\setminus G^\pard\mid C\in \Gamma\} 
\setminus \{\ecl\}$

$\okk G = G\setminus G^\pard = \bigcup\okk\Gamma$.
\\
$G^\pard$ contains all provably paradoxical atoms, while its
complement $\okk\gr$ could be taken as the atomic extension of the
consistency-operator, if we were aiming at a logic of formal
inconsistency. As we will see, it coincides with the domain of the maximal consistent subdiscourse, $\invE\gr[S]$, for any
$S\in Mod(\gr)$.
For now, we only note that
$\gr^\pard$ is $\E\gr$-closed, Fact \ref{AgClosed}, and that 
$\okk\Gamma$ remains consistent alongside $G^\pard$ and conservative
over $\Gamma$ with respect to the nonparadoxical atoms $\okk G$, as
made precise by Fact \ref{le:above}.
\begin{Fact}\label{AgClosed}For every $\Gamma$ and $x\in\gr:$ 
%with $\E\gr(x)=\{y_1...y_n\}:$ 
$\Gamma\pro \pard(x) \Impl \forall y\in \E\gr(x): \Gamma\pro\pard(y)$.
\end{Fact}

\begin{Proof}
 $\Gamma\pro x$ and the axiom $\nandp xy_i$, for each $y_i\in
 \E\gr(x)$, give $\Gamma\pro \ovr y_i$. Resolving then $xy_1...y_n$ with
 $\ovr x$ and $\ovr y_j$, for all $j\not=i$, gives $\Gamma\pro
 \ovr y_i$ for each $y_i\in \E\gr(x)$.
\end{Proof}
\begin{Fact}\label{le:above} 
For $\Gamma$ with $\okk G\not=\eset$:
\begin{enumerate}%\MyLPar
\item\label{it:a} $\forall  A\in\okk\Gamma:\Gamma\pro A$, so
$\forall C \subseteq \okk{\any G}: \okk\Gamma\pro C\Impl \Gamma\pro C$.

\item\label{it:parsem} $\okk\Gamma\not\pro\ecl$.

\item\label{it:single} $\forall x\in \okk G: \okk\Gamma\pro x \Iff
\Gamma\pro x$ and $\okk\Gamma\pro \ovr x \Iff \Gamma\pro \ovr x$.

\item\label{it:un}
$\exists x\in \okk G: \okk\Gamma\not\pro \ovr x$, hence also $\Gamma\not\pro\ovr x$.

\item\label{it:notint} $\forall x\in \okk G: \okk\Gamma\not\pro \ovr x
 \Impl \E\gr(x) \subseteq \okk G$ (when $\Gamma$ is a graph).
\end{enumerate}
\end{Fact}
\begin{Proof}
\ref{it:a}. Since for each atom $a\not\in \okk\gr$, both $\Gamma\pro a$ and
$\Gamma\pro\ovr a$, such atoms can be resloved away from every clause
of $\Gamma$. Each clause of $\okk\Gamma$ is obtained exactly
by such an operation.
\\[.5ex]
\ref{it:parsem}. $\okk\Gamma\pro\ecl \by\Impl{\ref{RES}.\ref{xnotx}} \exists x\in \okk\gr:\okk\Gamma\pro x\land \okk\Gamma\pro\ovr x
\by\Impl{\ref{it:a}.} \Gamma\pro x\land \Gamma\pro\ovr x \Impl x\not\in \okk\gr$.
\\[.5ex]
\ref{it:single}. Implications to the right follow by point \ref{it:a}, while to the left by Lemma \ref{lemma} and point \ref{it:parsem}.
\\[.5ex]
\ref{it:un}. If $\forall x\in \okk\gr: \okk\Gamma\pro \ovr x$, then also 
$\forall y\in \gr: \Gamma\pro \ovr y$, and then  
$\forall y\in \gr: \Gamma\pro y$, contradicting $\okk\gr\not=\eset$.
\\[.5ex]
\ref{it:notint}. If $x\in\okk\gr$ has a $y\in \E\gr(x)\cap \gr^\pard$, then $\Gamma\pro y\Impl \Gamma\pro \ovr x \by\Impl{\ref{it:single}.} \okk\Gamma\pro \ovr x$.
\end{Proof}
When $\Gamma$ represents a graph $\gr$, $\okk\Gamma$ is almost the
theory of its induced subgraph $\okk G$, except for a difference at
its border $brd(\okk G) = \{x\in \okk G\mid \E\gr(x) \not\subseteq
\okk G\}$.  For instance, for our discourse 
$\ourGr = \draw{a' \ar@<2pt>[r] & a\ar@<2pt>[l] & b\ar[l]\ar[r] & 
c\ar[r]& d\ar[r] & e\ar@/_/[ll]}:$ 
\[\begin{array}{lcl}
\Delta &=& \{aa', \nandp aa',bac,\nandp ba, \nandp bc, cd,\nandp cd, de, \nandp de, ec, \nandp ec\}
%, ac, \nandp ab, \nandp bc, \nandp ac, cde, \nandp cd, \nandp de, ef, \nandp ef\}
%\{ab, bc, ac, \nandp ab, \nandp bc, \nandp ac, cde, \nandp cd, \nandp de, ef, \nandp ef\}
\\
\ourGr^\pard &=& \{c,d,e\} \\
\okk\ourGr &=& \{a',a,b\}\hfill \mbox{-- or the induced subgraph }\draw{a'\ar@<2pt>[r]& a\ar@<2pt>[l] & b\ar[l]}
\\
\okk\Delta &=&  \{aa', \nandp aa',ba,\nandp ba, \ovr b\}
%\{de, \nandp de, ef, \nandp ef\} \cup \{\ovr d\}
\\
brd(\okk\ourGr) &=& \{b\}
\\
\cthr(\okk\ourGr) &=&  \{aa', \nandp aa', ba, \nandp ba\} 
\mbox{\ -- clausal theory of the induced  subgraph $\okk\ourGr$}
\end{array}
\] As in this example so generally, border vertices enter as negative
literals into \(\okk\Gamma = \cthr(\okk\gr)\cup (brd(\okk G))^-\),
according to Fact \ref{le:above}.\refp{it:notint}. One can thus view
$\okk\Gamma$ as (the theory of) the subgraph induced by $\okk G$, with
a loop added at each border vertex. It is consistent, Fact
\ref{le:above}.\refp{it:parsem}, so its models are kernels of
$\okk{\gr}$ excluding border vertices:
\eq{CMod(\okk\Gamma) = \{L\in \Ker{\okk{\gr}}\mid brd(\okk G)\subseteq
\invE\gr(L)\}.\label{CModGokk}}
These classical models of $\okk\Gamma$ are
 actually the models $Mod(\Gamma)$ from 
 \refp{Mod}. %for the graph $\gr$ of $\Gamma$. 
To show this, we first register a preliminary observation.
\begin{Fact}\label{semi}
$CMod(\okk\Gamma)\subseteq SK(\gr)$.
\end{Fact}
\begin{Proof}
Since each $L\in CMod(\okk\Gamma)$ is a kernel of $\okk\gr$, it is
obviously independent in $\okk\gr$ and, since $\okk\gr$ is an induced
subgraph of $\gr$ (only with additional loops), so $L$ is independent
also in $\gr$. Hence $\invE\gr(L)\subseteq \gr\setminus L$. 

By soundness of \RES, if $a\in L$ then $\okk\Gamma\not\pro \ovr a$, 
which by Fact \ref{le:above}.\refp{it:notint} means that $\E\gr(a)\subseteq \okk\gr$, so that 
$\E\gr(a) = \E{\okk\gr}(a) \subseteq \invE{\okk\gr}(L)$, where the inclusion follows since $L\in Ker(\okk\gr)$. Obviously, $\invE{\okk\gr}(L)\subseteq  \invE{\gr}(L)$. Combining the two yields
$\E\gr(L)\subseteq \invE\gr(L)\subseteq \gr\setminus L$, i.e., $L\in SK(\gr)$. 
\end{Proof}
%Each $S\in SK(\gr)$ determines a 3-partition $\alpha_S =
%  \<S,\invE\gr(S),\gr \setminus \invE\gr[S]\>$, so writing
%  $S\models\Gamma$ or $S\in pSK(\gr)$, we mean $\alpha_S\models\Gamma$
%  or $\alpha_S\in pSK(\gr)$.
Writing $S\models\Gamma$ or $S\in pSK(\gr)$, for an $S\in SK(\gr)$,
  we mean $\alpha_S\models\Gamma$ or $\alpha_S\in pSK(\gr)$, for $\alpha_S =
  \<S,\invE\gr(S),\gr \setminus \invE\gr[S]\>$.
The bijection in the theorem, $CMod(\okk\Gamma) \simeq Mod(\gr)$, means then $\invE{\okk\gr}(S)=\invE\gr(S)$, for
  relevant $S$, and
\begin{itemize}
\item[($\supsetsim$)] $\forall \<S,\invE\gr(S),\gr\setminus \invE\gr[S]\>\in Mod(\gr):
\<S,\invE{\gr}(S)\>\in CMod(\okk\gr)$ and  \vspace*{.5ex}
\item[($\subsetsim$)]
 $\forall \<S,\invE{\okk\gr}(S)\> \in CMod(\okk\Gamma):
\<S,\invE\gr(S),\gr\setminus \invE\gr[S]\>\in Mod(\gr)$. 
\end{itemize}
\begin{Theorem}\label{not-main}
For every theory $\Gamma$ with graph $\gr: CMod(\okk\Gamma) \simeq Mod(\gr)$. 
\end{Theorem}
\begin{Proof}
{\bf{($\supsetsim$).}} follows using Fact 
  \ref{sound} but first we show that (1a) if $S\in Mod(\gr)$ then
  $S\models\Gamma$.  This and (1b) hold actually for every
  $S\in pSK(\gr)$:
  \\
  {\bf (1a)} $S\in pSK(\gr) \Impl S\models\Gamma$.  For each
  $y\in \E\gr(x)$, i.e., $\ovr{xy}$, we have one of four cases, each
  yielding $S\models \nandp xy$:
\begin{itemize}
\item[(i)] $y\in S \Impl x\in \invE\gr(S)$, 
\item[(ii)] $x\in S \Impl y\not\in S$ and, since $S\in SK(\gr)$, $y\in \invE\gr(S)$, 
\item[(iii)] $x\in \invE\gr(S)$ or $y\in \invE\gr(S)$,   
\item[(iv)] $\{x,y\}\subseteq \gr\setminus \invE\gr[S]$.
\end{itemize}
For each $\E\gr[x] = \{x\}\cup Y$, we have one of five cases, each
giving $S\models\E\gr[x]$:
\begin{itemize}
\item[(i)] $\exists y\in Y: y\in S$, 
\item[(ii)] $x\in S$, 
\item[(iii)] $x\in \invE\gr(S) \impl \exists y\in Y: y\in S$, 
\item[(iv)] $\exists y\in Y: y\in \invE\gr(S)$, since $S\in pSK(\gr)$, so $x\in \invE\gr[S]$ and 
 $S\models \E\gr[x]$ by (ii) or (iii), 
\item[(v)] $\E\gr[x]\subseteq \gr\setminus \invE\gr[S]$. 
\end{itemize}
{\bf (1b)} $\forall S\in pSK(\gr): \Gamma\pro\pard(x) \Impl x\not\in \invE\gr[S]$. 

By (1a) $S\in pSK(\gr) \Impl S\models\Gamma$ so, by soundness Fact
\ref{sound}, for any clause $C:\Gamma\pro C\Impl S\models C$. Hence,
$\Gamma\pro\pard(x) \Impl S\models \pard(x)$, i.e.,
$x\in \gr\setminus \invE\gr[S]$.
\\
{\bf (1c)} We show that $S\in Mod(\gr)$ satisfies \refp{CModGokk}.
By Fact \ref{le:above}.\refp{it:parsem},
$\okk\Gamma\not\pro\ecl$ so, by Fact \ref{RES}.\refp{ref}, there is
some $K\in \sols{\okk\gr}$, i.e., one with $\invE\gr[K] = \okk\gr$.
By (2a-b) below, $K\in pSK(\gr)$; since $S\in mpSK(\gr)$ (as
$S\in Mod(\gr)$ and \refp{mpSK}), so
$\invE\gr[S]\supseteq \invE\gr[K]=\okk\gr$.  By (1b),
$x\in \invE\gr[S]\Impl x\not\in \gr^\pard\Impl x\in \okk\gr$, i.e.,
$\invE\gr[S]\subseteq \okk\gr$, so that $\invE\gr[S] = \okk\gr$.  By
\refp{ref:closure}, $\invE\gr(\okk\gr)\subseteq\okk\gr$, so since
$\okk\gr$ is induced subgraph of $\gr:\invE{\okk\gr}(S) = \invE\gr(S)$
and then, since $S\cap \invE\gr(S)=\eset$ as $S\in SK(\gr)$, so
$\invE{\okk\gr}(S) = \okk\gr\setminus S$ -- showing that
$S\in \sols{\okk\gr}$.

If $x\in brd(\okk\gr)$ then $\Gamma\pro \ovr x$, so $S\models\ovr x$,
i.e., $x\not\in S$ and, since
$S\in \sols{\okk\gr}: x\in \okk\gr\setminus S = \invE\gr[S]\setminus S
= \invE\gr(S)$. Thus $brd(\okk\gr)\subseteq \invE\gr(S)$.
\\[1ex]
{\bf ($\subsetsim$).}  {\bf (2a)} $CMod(\okk\Gamma)\subseteq SK(\gr)$ is Fact \ref{semi}.
%is Fact 5.6 from \cite{RIP}.
\\[1ex]
{\bf (2b)} By Fact \ref{AgClosed}, 
$\Gamma\pro \pard(x) \Impl \forall y\in \E\gr(x): \Gamma\pro\pard(y)$,
so $\E\gr(\gr^\pard)\subseteq \gr^\pard$. Hence
$\invE\gr(\gr\setminus \gr^\pard)\subseteq (\gr\setminus \gr^\pard)$,
i.e., $\invE\gr(\okk\gr)\subseteq \okk\gr$ and so
\\
$\begin{array}{rcl}
\invE\gr(\invE\gr[\okk\gr]) & = & \invE\gr(\invE\gr(\okk\gr) \cup
\okk\gr) 
\\ & = & \invE\gr(\invE\gr(\okk\gr)) \cup \invE\gr(\okk\gr)
\subseteq \invE\gr(\okk\gr) \cup \okk\gr = \invE\gr[\okk\gr]. \end{array}$
\\[1ex]
{\bf (2c)} When $S\in \sols{\okk\gr}$ then $\invE\gr[S]=\okk\gr$ and
$S\in pSK(\gr)$, by (2a-b).  If $S\not\in mpSK(\gr)$, i.e.,
$\exists R\in pSK(\gr): \invE\gr[R]\not\subseteq \invE\gr[S]$, Lemma
\ref{lemma:main0} yields a strict extension
$\invE\gr[Q]\supset \invE\gr[S]$. This requires adding some 
$E\subseteq \gr^\pard$, but by (1b) no $e\in E$ can belong to any
$\invE\gr[Q]$ with $Q\in pSK(\gr)$, since $\Gamma\pro \pard(e)$.
\end{Proof}
The theorem implies that
 $\<\alpha^\true,\alpha^\false,\alpha^\pard\>\in Mod(\gr)$ iff
 $\<\alpha^\true,\alpha^\false\>\in CMod(\okk\Gamma)$, giving the
 middle equality: $\gr\setminus \alpha^\pard =
 \alpha^\true\cup\alpha^\false = \okk\gr = \gr\setminus \gr^\pard$,
 which yields: 
\eq{\forall \alpha\in Mod(\gr): \alpha^\pard = \gr^\pard.\label{semsyntbot} } 
Thus \RES\ proves both $a$
 and $\ovr a$ exactly for the atoms $a$ falling outside the healthy, boolean
 domain of every model.  The following completeness of \RES\ for
 $\models$ is the counterpart of its classical completeness for
 $\models_c$ from Fact \ref{RES}.\refp{wCompl}.
\begin{Corollary}\label{cor:main}
For every $\Gamma$ and clause $A\subseteq \any \gr:$

\(
\begin{array}{ll}
\Gamma \models  A\ \  \Ifff\ \ 
A \subseteq G^\pard \not = \eset\ \text{ or }\  %\\ 
\exists B\not=\ecl: B \subseteq A \cap \okk{\any \gr} \land \Gamma \pro  B. 
\end{array} 
\) 
\end{Corollary} \vspace*{-1ex}

\begin{Proof} 
  $\Rightarrow$) Assume that $\Gamma \models A$. If $A \subseteq
  \any G^\pard$ then $G^\pard \not = \eset$ by (\ref{parasem}) and we are
  done. If $A \not \subseteq \any G^\pard$, i.e., $C = A \setminus \any G^\pard
  \not = \eset$, then $C \subseteq A \cap \okk{\any \gr}$ and, by (\ref{parasem}), 
  $\Gamma \models C$ (since $\forall\alpha \in Mod(\Gamma):
  G^\pard = \alpha^\pard$ by \refp{semsyntbot}). Hence, by Theorem
  \ref{not-main}, $\okk \Gamma \models_c C$ which, by Fact
  \ref{RES}.\refp{wCompl}, implies $\exists B\subseteq C :
  \okk\Gamma\pro B$. By Fact \ref{le:above}.\refp{it:parsem},
  $B\not=\ecl$, while by \ref{le:above}.\refp{it:a}, $\Gamma\pro
  B$, yielding the conclusion, since $B \subseteq A \cap \okk{\any \gr}$. 
\\[1ex] 
  $\Leftarrow$) If $A \subseteq \gr^\pard \not = \eset$ then $\Gamma
  \models A$ directly by (\ref{parasem}). Otherwise, as $\Gamma \pro
  B$ so $\Gamma \models B$ by Fact \ref{sound}. Since $\eset\not=B
  \subseteq A\cap \okk{\any \gr}$, so $\Gamma \models A$ by (\ref{parasem}).
\end{Proof}
Comparison with Fact \ref{RES}.\refp{wCompl}
 shows that while in
classical logic $\ecl$ witnesses to every clause by {\it Ex Falso}, in
our logic {\it all} consequences of a theory have nonempty provable
witnesses: $\Gamma\models\ecl$ iff $\gr^\pard\not=\eset$, i.e.,
$\exists x: \Gamma\pro x$ and $\Gamma\pro\ovr x$, while
$\Gamma\models A\not=\ecl$ is witnessed only by a proof from $\Gamma$
of either $\pard(a)$, for all $a\in A$, or of some nonempty subclause
$B\subseteq A$.

The implication to the left in Corollary \ref{cor:main} specifies more
closely the patterns of weakening admisible in our logic. Not only {\it
  Ex Falso} is excluded, but weakening has to preserve, so to say, the
reason of satisfaction of its premise. The first case (right-to-left
implication from the first disjunct in Corollary \ref{cor:main}) allows to
weaken the empty clause only to a disjunction of literals involved in
inconsistency, reflecting the irrelevance of the empty clause for the
consistent subdiscourse. The second case ensures that the resulting
$A$ contains a healthy (i.e., not involving any paradox)
$B\subseteq A$ witnessing to its satisfaction. Thus, adding to \RES\
either both (aW) and (bW) or only (cW)
\[
(aW)\ \prr{\Gamma\pro  B}{\Gamma\pro  B \cup 
  C}B\not\subseteq\any\gr^\pard \hspace*{1em} (bW)\
\prr{\Gamma\pro \pard(a)\ \ \forall a\in A\not=\eset}{\Gamma\pro 
  A}\hspace*{2em} (cW)\ \prr{\Gamma\pro B}{\Gamma\pro B \cup C} 
\]
yields a sound and strongly complete system for classical logic, with (cW), or
 for our paraconsistent logic, with (aW) and (bW). For consistent
 $\Gamma$, the two coincide. Since
 $\gr^\pard=\eset$, (bW) is inapplicable, while (aW) becomes exactly
 (cW), since then every $B\not\subseteq \gr^\pard$.
% the case of {\it Ex Falso}, i.e., provability of  $B=\ecl$, does not occur.

The conditions of (aW) and (bW), reflecting the right 
side of Corollary \ref{cor:main}, prevent uncontrolled mixing of
consistent and inconsistent elements. As an example, consider the usual 
 derivation of Lewis' ``paradox'', recasting {\it Ex Falso} using
disjunctive syllogism. Assuming $\Gamma\pro a$ and $\Gamma\pro\ovr a$,
we have

\prr{\prr{\Gamma\pro a}{\Gamma\pro ab}\ (cW) \ \ \prline{\Gamma\pro
    \ovr a}}{\Gamma\pro b}\ (DS)
\\
The step (DS) is an instance of resolution, but since
$\{a\}\subseteq \gr^\pard$, (aW) can not be applied instead of
(cW).

\section{A remark on relevance}\label{rel}
Rule (aW) is not very effective, requiring to find $b\in B$ for which
$b$ or $\ovr b$ is not provable, while rule (bW) joins nonsensical
premises into a nonsensical disjunction. These rules, extending \RES\ to a
strongly complete system, provoke also the question: what do we lose by 
dropping them? \RES\ remains complete in the sense of Corollary
\ref{cor:main}; if $a\lor b$ is a consequence of our theory, we may be
able to prove $a$, but not necessarily $a\lor b$. Is it such a loss? If we know that
$a$ is true, it is not particularly enlightening that so is $a\lor
b$. In this way, the tautology $\neg a\lor a$ can be diluted to, say, 
\[(i)\ \neg a \lor a \lor \neg b.\] With implication $x\impl y$
defined as $\neg x\lor y$, (i) represents some fallacies of relevance,
typical for material implication, e.g., $a\impl (b\impl a)$ and $\neg
a \impl (a\impl \neg b)$, for strict implication, e.g., $b \impl
(a\impl a)$ or $b\impl (\neg a\lor a)$, or for the intutionistic one,
$\neg(a\impl a) \impl b$. \RES\ is immune against such fallacies,
which arise in the language of clauses from disjunctive weakening.
Note that while this rule may be indispensable in various systems of
Gentzen or natural deduction, in \RES\ it is not, since here its only
contribution is dilution of provable facts. Dropping all forms of
weakening, we lose such diluted consequences, but gain
relevance, simplifying reasoning at the same time. The following 
remarks  elaborate these gains. Without aiming at a logic of
relevance, they only identify its elements in the present setting.

Let us first note that the language of clauses, viewed as sets of
literals, identifies, for instance, $a\impl (a\impl b)$ and $a\impl
b$, representing both as $\ovr ab$ and enforcing their
equivalence. One could see it as an unfortunate equivocation but, with
the present semantics, it only eliminates spurious
syntax.\footnote{\RES\ with this language can be seen as a restriction
of $RMI_{\scriptscriptstyle{\stackrel\sim{+}}}$ with the usual
sequential syntax, from \cite{AvronResol}. The equivalence of
$RMI_{\scriptscriptstyle{\stackrel\sim{+}}}$ and
$RMI_{\scriptscriptstyle{\stackrel\sim{\impl}}}$ shows that the
former's disjunction (specializing to ours in the present context)
yields in the latter a satisfactory relevant implication, defined by
$A\impl B = \neg A \lor B$. To handle subtler aspects of relevant
implication, such an extension of the syntax and \RES, along with the
associated semantic adjustments, might be needed. }

The absence of weakening prevents \RES\ from deriving, besides
fallacies exemplified above, also some other 
problematic implications, for instance, $(a\impl b)\land (c\impl
d)\not\pro (a\impl d) \lor (c\impl b)$. From the assumption, represented by the 
clauses $\ovr ab, \ovr cd$, nothing follows by resolution,
in particular, not the undesired conclusion, $\ovr ad\ovr cb$, which
only weakens either premise.

\RES\ enjoys a specifically relevant form of deduction theorem.  If
$\Gamma, x\pro y$, it does not follow that $\Gamma\pro x\impl y$,
i.e., $\Gamma\pro \ovr xy$, without further ado. It may namely happen
that $\Gamma\pro y$ without using $x$, and $\Gamma\not\pro \ovr
xy$. By Fact \ref{RES}.\refp{aux}, however, if the proof $\Gamma,x\pro
y$ requires $x$, that is, $\Gamma\not\pro y$, then indeed $\Gamma\pro
\ovr xy$.

Just as our logic is concerned with consequences of a given
theory rather than with tautologies, relevance is judged relatively to the
actual context.  From some $\Gamma$, (i) may be provable and from
others it may not be. For instance

(i1) $\Gamma_1 \pro ca\ovr b$ and $\Gamma_1\pro \ovr a a\ovr b$, where 
 $G_1= ...c \leftarrow a \into d\into b...$, while 

(i2) $\Gamma_2 \pro \ovr a a$\ \ but\hspace*{.12em}  $\Gamma_2\not\pro \ovr a a \ovr b$, where  
 $G_2= ...c \leftarrow a \hspace*{1.5em} d\into b...$,
\\
where ``$...$'' indicates a sur\-round\-ing graph, in which $a$ and
$d$ have no incident edges except the indicated ones.  The graph
syntax makes explicit the reference structure, allowing to read an
edge $x\into y$ as $x$ ``referring to'' $y$, by negating
$y$. Reference involves relevance, so that $y$, negated by $x$, is
relevant for $x$.  The atoms reachable by arbitrary paths from a given
$x$ are then indirect references, relevant for $x$. This form of
relevance is thus transitive, but context dependence and the graph
syntax justify this fact. In (i1), $b$ is relevant for $a$, because of
their connection through $d$. Making $b=\false$, forces
$a=\false$. This notion, involving only context dependent interaction between
some truth-values of $b$ and $a$, deviates from those attempting to
capture some general meaning-connections, for instance, by
variable-sharing. The two facts in ``Water freezes or temperature
is above 0$^\circ$C'' do not share any variable but are connected
by the background knowledge. The meaning-connections, already between
atomic facts, are defined by the context $\Gamma$ and reflected by its
undiluted consequences.  Provability of $\ovr aa\ovr b$ in (i1)
suggests such a connection of $b$ and $a$, while its unprovability in
(i2) witnesses to its absence.

Relevance is more than reference -- typically, it is symmetric. In
our case, this amounts to following the paths also in the direction opposite
to the edges:\vspace*{1ex}

\(\begin{array}{l@{\hspace*{3em}}l}
\mbox{(a) This and the next statmenet are false.} &  a \iff  \neg a \land \neg b \\[1.5ex]
\mbox{(b) The next statmenet is false.}           &  b \iff  \neg c \\[1.5ex]
\mbox{(c) The previous statmenet is false.}       &  c \iff  \neg b
\end{array} \hspace*{2em}
\raisebox{1cm}{\xymatrix@R=.4cm{a \ar@(ru,rd) \ar[d] \\ b \ar@/_/[d] \\ c \ar@/_/[u]} }
\)
\\
In the only model, $b=\true$ and $a=c=\false$. But if we remove
the loop at $a$, there will be also another model with $a=c=\true$ and
$b=\false$.  One could say that the loop at $a$ forces $b=\true$,
since otherwise inconsistency would arise. Not only $b$, to which $a$
refers, is relevant for $a$, but also $a$ is relevant for $b$, forcing
it to be true.  In this sense, all ``connections'' -- paths
in the underlying graph (the undirected graph, obtained by
forgetting directions of the edges) -- mark potential relevance. It is
thus both symmetric and transitive, so its horizon for an atom $x$
is the strong component of the underlying graph, containing $x$. 
Indeed, \RES\ proves some clause with literals over atoms $x$
 and $y$ if and only if $x$ and $y$ belong to one such component.\footnote{Distinct strong components have disjoint
 alphabets, with no axiomatic clauses containing literals from both;
 hence no provable clause can contain literals from both. Conversely,
 if there is a path $\<x_1,x_2,...,x_n\>$ in the underlying graph, then
 for each pair $x_i,x_{i+1}$, $1\leq i<n$, there is a negtive clause
 $N_i=\ovr{x_i}\ovr x_{i+1}$ and a positive one
 $P_i=...x_ix_{i+1}...$. Resolving $P_1$ with $N_2$ yields a clause
 $C_3$ containing $x_1$ and $\ovr x_3$. Resolving $C_3$ with $P_3$
 gives $C_4$ containing $x_1$ and $x_4$. A clause $C_n$, obtained after $n-2$ steps, 
 contains $x_1$ and either $x_n$ or $\ovr x_n$.}

Still, this notion of relevance %obtained by requiring containment in a strong component 
is too weak, as exemplified by the clause $\ovr aa\ovr
b$ from (i1). Signaling a connection between $a$ and $b$, it does not
impose any dependencies between the truth-values of $a$, $b$ and the
whole clause.  Its proof in \RES\ is not a mere dilution of $\ovr aa$,
but the clause itself is such a dilution, which is always true because
so is its subclause $\ovr aa$.  The sense of relevance should preclude
us from saying such things which, even if true, do not add anything to
something which is already said more concisely.

Let us therefore consider a clause $C$ {relevant}, in a given context
$\Gamma$, denoted $\Gamma\models_r C$, if $\Gamma\models C$ and for
each nonempty $B\subset C: \Gamma\not\models B$. For unhealthy atoms,
such relevant clauses are the units $a$ and $\ovr a$, and no other relevant
clause contains such atoms. This reminds, once more, of the
irrelevance of meaningless/inconsistent elements for the healthy part
of a discourse.  For each 2-partition $C=AB$ of a nonunit clause $C$
(with at least two literals), we can write $\Gamma\models C$ as
$\Gamma\models \band A^- \impl \bor B$. If $\Gamma\models_r C$, this
implication is not satisfied, so to say, vacuously, by $\Gamma$ having
no models satisfying $\band A^-$, but by $\Gamma$ having such models,
all satisfying also $\bor B$. A relevant $C$ witnesses thus to the
influence, which some truth-values of any 
%relevance of every 
nonempty proper subclause $A\subset C$ %for 
have for 
its complement
$B=C\setminus A$: whenever $\band A^-=\true$, then $\bor B=\true$, and
there are cases when $\band A^-=\true$. This last proviso is the
element of relevance.  The relation is obviously symmetric: if
$\band A^-\impl \bor B$, then also $\band B^-\impl \bor A$.
%$A$ and $B$ are relevant if $\Gamma\models_r AB$.

This is a much stronger notion than mere membership in the same
component of the underlying graph. In the example (i1),
$\Gamma_1\not\models_r\ovr aa\ovr b$ since
$\ovr aa\subseteq \ovr a a\ovr b$, so the fact that
$\Gamma_1\models\ovr aa\ovr b$ does not imply relevance in this strong
sense. On the other hand, $\Gamma_1\pro \ovr ab$ and 
$\Gamma_1\models_r\ovr ab$, so $\ovr a$ and $b$ are relevant for
each other: $\Gamma_1$ has then models with $b=\false$, in which
$a=\false$, and models with $a=\true$, in which $b=\true$. If
$\Gamma_1\not\models_r \ovr ab$ then $\ovr a$ and $b$ would cease to be
mutually relevant, since then either $\Gamma_1\models\ovr a$ or
$\Gamma\models_r b$, irrespectively of each other.

This relevance relation is not transitive. In 
\draw{\gr= ...d &  a\ar[l]\ar[r] & b\ar[r] \ar@/^/[rr] & c... & e...} 
%\draw{...a\ar[r] & b\ar[r] \ar@/^/[rr] & c... & e...} 
\\
(with appropriate assumptions about the undisplayed part,
e.g., all $...$ starting mutually disjoint infinite paths), $b$ is
relevant for $a$, since $\Gamma\models_r \nandp ab$, and $c$ is
relevant for $b$ since $\Gamma\models_r\nandp bc$. But $c$ is no
longer relevant for $a$, because any value at $c$ can be accompanied
by any value at $a$. However, $c\lor e$ is re\-le\-vant for $a$, since
$\Gamma\models_r \ovr ace$.%\noo{ ($\Gamma$ does not satisfy any binary clause over $a$ and $c$.)} 

The relevant clauses are exactly  the 
 minimal provable ones:~\footnote{First, for each nonunit clause $C\in
Min(\Gamma)$ and nonempty $B\subset C: \okk\Gamma,B^-\not\pro\ecl$.
This follows because, by Fact \ref{le:above}.\refp{it:parsem},
$\okk\Gamma\not\pro\ecl$ and then, by Fact \ref{RES}.\refp{aux},
$\okk\Gamma,B^-\not\pro\ecl$ iff $\forall E\subseteq
B:\okk\Gamma\not\pro E$, which holds since $C\in Min(\Gamma)$. The
main claim is trivial for unit clauses.  Assuming now a nonunit
$AB=C\in Min(\Gamma)$, the observation above implies existence of a
classical model of $\okk\Gamma$ (i.e., a model of $\Gamma$),
satisfying $B^-$, so that $\Gamma\not\models B$. Since $\okk\Gamma\pro
C$, so $\Gamma\pro C$ by Fact \ref{le:above}.\refp{it:a}, and hence
$\Gamma\models C$ by soundness. So $C\in Min(\Gamma)\Impl
\Gamma\models_r C$. Conversely, if for a nonunit $C:
\Gamma\models_rC$, then $C\cap \any\gr^\pard=\eset$, so that
$\okk\Gamma\models_c C$. Then $\okk\Gamma\pro B$ for some $B\subseteq
C$ by Fact \ref{le:above}.\refp{wCompl}. But since $\Gamma\not\models
B$, for every $B\subset C$ so, by soudness, $\Gamma\not\pro B$, and
hence $\okk\Gamma\not\pro B$ by \ref{le:above}.\refp{it:a}, so
$\okk\Gamma\pro C$ and $C\in Min(\Gamma)$.}

\(\begin{array}{rcl}Min(\Gamma) &=& \{C\mid \Gamma\pro C \land \forall
B\subset C: B\not=\ecl \Impl \Gamma\not\pro B\}
\\
& = & \{C\mid \okk\Gamma\pro C \land \forall
B\subset C: \okk\Gamma\not\pro B\} \cup %\bigcup_{\Gamma\pro\pard(a)}\{a,\ovr a\}.
\bigcup\{\{a,\ovr a\}\mid \Gamma\pro \pard(a)\}.
\end{array}\)
\\
Thus each \RES\ provable clause $C$ represents a form of
relevance. Even if $C\not\in Min(\Gamma)$, all its atoms belong to one
component of the underlying graph and, moreover, there is some
nonempty $D\subset C$ with $D\in Min(\Gamma)$, witnessing to the
determination -- in some occurring cases -- of the truth-value of one
part of each 2-partition of $D$ by the truth-value of the other.

\bibliographystyle{plain} %{jflnat} 
\bibliography{../biblio}

\end{document}